\input amstex
\documentstyle{amsppt}
\magnification=\magstep1

\baselineskip=13pt
\parskip3pt

\define\m1{^{-1}}
\define\ov1{\overline}
\def\gp#1{\langle#1\rangle}
\def\ul2#1{\underline{\underline{#1}}}

\def\Y{\text{ {\rm Y} }}
\def \iiitem{\par\indent\indent\hangindent 3\parindent\textindent}

\topmatter
\title
Structure of  normal twisted group rings
\endtitle
\author
VICTOR BOVDI
\endauthor
\address
       Victor Bovdi\newline
       Institute of Mathematics and Informatics\newline
       Lajos Kossuth University \newline
       H-4010  Debrecen, P.O.Box 12\newline
       Hungary
\endaddress
\email
vbovdi{\@}math.klte.hu
\endemail
\abstract
Let $K_{\lambda}G$ be the twisted group ring of a group $G$ over a
commutative ring $K$ with $1$, and let $\lambda$ be a factor set
($2$-cocycle) of $G$ over $K$. Suppose $f: G\to U(K)$ is  a map
from $G$ onto the group of units $U(K)$ of the ring $K$ satisfying
$f(1)=1$. If $x=\sum_{g\in G}\alpha_gu_g \in K_\lambda G$ then we
denote $\sum_{g\in G}\alpha_gf(g)u_g^{-1}$ by $x^f$ and assume that
the map $x\to x^f$ is an involution of $K_{\lambda}G$. In this
paper we  describe those groups $G$ and commutative rings $K$ for which
$K_{\lambda}G$ is  $f$-normal, i.e\.  $xx^f=x^fx$ for
all $x\in K_{\lambda}G$.
\endabstract
\thanks
Research supported by OTKA  No\.~T16432.
\endthanks

\endtopmatter

\document

\head {1. Introduction} \endhead

Let $G$ be a group and $K$ a commutative ring with unity.
Suppose that the  elements of the set\quad 
$
\Lambda= \{\lambda_{a, b} \in U(K)\mid a, b\in G\}$\quad  satisfy the condition
$$
\lambda_{a, b}\lambda_{ab, c}=\lambda_{b, c}\lambda_{a, bc} \tag 1 
$$
for all $a, b, c\in G$. Then $\Lambda$ will be called a {\it factor
system\/} ($2$-cocycle) of the group $G$ over the ring $K$.
The twisted group ring $K_{\lambda}G$ of $G$ over the commutative
ring $K$ is an associative $K$-algebra with basis $\{u_g \mid g\in
G\}$ and with multiplication defined distributively by
$u_gu_h=\lambda_{g, h}u_{gh}, $ where $g, h\in G$ and
$$
\lambda_{g, h}\in \Lambda=
\{\lambda_{a, b} \in U(K)\mid a, b\in G\}.
$$
Note that if $\lambda_{g, h}=1$ for all $g, h\in G$,
then $K_{\lambda}G\cong KG$, where $KG$ is the group ring of the
group $G$ over the ring $K$.

Properties of twisted group algebras and their groups of units
were studided by many authors, see, for instance, the paper by
S. V.~{\smc Mihovski} and J. M.~{\smc Dimitrova} \cite{1}.
Our aim is to describe the structure of   $f$-normal
twisted group rings. This result for  group rings was obtained in
\cite{2, 3}.

We shall refer to two twisted group rings  $K_{\lambda}G$ and
$K_{\mu}G$ as being diagonally equivalent if there exists a map
$\theta:G\to U(K)$   such that
$$
\lambda_{a, b}=\theta(a)\theta(b)\mu_{a, b}(\theta(ab))\m1.
$$
We
say that a factor system $\Lambda$ is normalized if it satisfies
the condition
$$
\lambda_{a, 1}=\lambda_{1, b}=\lambda_{1, 1}=1
$$
for all $a, b\in G$.

Hence, given $K_{\mu}G$ there always exists a diagonally equivalent
twisted group ring $K_\lambda G$ with factor system $\Lambda$
defined by $\lambda_{a, b}=\mu_{1, 1}\m1\mu_{a, b}$ such that
$\Lambda$ is normalized. From now on, all the factor systems
considered are supposed to be normalized.

The map $\phi$ from the ring $K_\lambda G$ onto $K_\lambda G$ is called
{\it an involution\/}, if it satisfies the conditions

\rm{(i)}    $\phi(a+b)=\phi(a)+\phi(b)$;\hskip 10pt
\rm{(ii)}   $\phi(ab)=\phi(b)\phi(a)$;\hskip 10pt
\rm{(iii)}  $\phi^2(a)=a$ \newline
for all $a, b\in K_\lambda G$.

Let $f:G\to U(K)$ be a map from the group  $G$ onto the group of
units $U(K)$ of the  commutative ring $K$, satisfying $f(1)=1$. For an
element $x=\sum_{g\in G}\alpha_gu_g \in K_\lambda G$ we define
$x^f=~\sum_{g\in G}\alpha_gf(g)u_g\m1 \in K_\lambda G$.

Let $x\to x^f$ be an involution of the twisted group ring
$K_\lambda G$. The twisted group ring $K_\lambda G$ is called
{\it $f$-normal\/} if
$$
xx^f=x^fx \tag2
$$
for all $x\in K_\lambda G$.

Recall that a $p$-group is called {\it extraspecial\/}
(see \cite{4}, Definition III.13.1)
if its centre, commutator subgroup and Frattini subgroup are equal and
have order $p$.

\proclaim {Theorem}
Let $x\to x^f$ be an involution of the twisted group ring
$K_\lambda G$.  If the ring $K_\lambda G$ is $f$-normal then
the group $G$ and the ring $K$ satisfy one of the following
conditions:

\smallskip
\item {\rm 1)}
$G$ is  abelian  and the factor system is symmetric, i.e\.
$\lambda_{a, b}=\lambda_{b, a}$ for all $a, b \in G$;

\smallskip
\item {\rm 2)}
$G$ is an abelian group of exponent $2$ and the factor
system  satisfies
$$
(\lambda_{a, b}-\lambda_{b, a})(1+f(b)\lambda_{b, b}\m1)=0 \tag 3 
$$
for all $a, b\in G$;

\smallskip
\item {\rm 3)}
$G=H\rtimes C_2$ is a semidirect product of an
abelian group $H$ of exponent not equal to $2$ and  $C_2=\gp{a\mid a^2=1}$
with $h^a=h\m1$ for all $h\in H$, the factor system of $H$ is
symmetric, $f(a)=-\lambda_{a, a}$ and
$$
\lambda_{a, h}=f(h)\lambda_{h, h\m1}\m1\lambda_{h\m1, a}, \quad
\lambda_{h, a}=f(h)\lambda_{h, h\m1}\m1\lambda_{a, h\m1};\tag 4 
$$

\smallskip
\item {\rm 4)}
$G$ is a hamiltonian $2$-group and the factor system  satisfies
\iiitem {\rm 4.i)}
for all noncommuting $a, b \in G$
$$
\lambda_{a, b}=f(a)\lambda_{a, a\m1}\m1\lambda_{b, a\m1}=
f(b)\lambda_{b, b\m1}\m1\lambda_{b\m1, a};\tag 5 $$

\iiitem {\rm 4.ii)}
$\lambda_{g, h}=\lambda_{h,g}$ for any $h\in C_G(\gp{g})$ and
$f(c)=\lambda_{c, c}$ for every $c$ of order $2$;

\smallskip
\item {\rm 5)}
$G=\Gamma \Y C_4$ is a central product of a hamiltonian $2$-group $\Gamma$
and a cyclic group $C_4=\gp{d\mid d^4=1}$  with $\Gamma'=\gp{d^2}$.
The factor system satisfies {\rm (5)} and
$$
\lambda_{b, a}\lambda_{ba, d}+f(d)\lambda_{d, d\m1}\m1\lambda_{a, b}
\lambda_{ab, d\m1}=0, \tag 6 $$

\item {}
where $a,b\in \Gamma$, $a^4=b^4=1$ and $[a,b]\ne 1$;

\smallskip
\item {\rm 6)}
$G$ is either $E\times W$ or $(E \Y C_4)\times W$, where $E$ is an
extraspecial 2\hbox{-}group, $E \Y C_4$ is  the central product of $E$  and
$C_4=\gp{c\mid c^4=1}$  with
$E'=\gp{c^2}$ and  $\exp(W)|2$. The factor system satisfies:
\iiitem {\rm 6.i)}
If $a\in G$ has order $4$ then $\lambda_{a,h}=\lambda_{h, a}$ for all
$h\in C_G(\gp{a})$;
\iiitem {\rm 6.ii)}
if $\gp{a, b}$ is a quaternion subgroup of
order $8$ of $G$ then the properties {\rm (5)} and {\rm (6)} are satisfied for
every $d\in C_G(\gp{a, b})$ of order $4$,
and $f(v)=\lambda_{v, v}$ for all $v\in C_G(\gp{a, b})$ of order $2$;
\iiitem {\rm 6.iii)}
if  $\gp{a, b \mid a^4=b^2=1}$ is the dihedral group of order $8$, then
$f(b)=-\lambda_{b, b}$  and the properties {\rm (4)}, {\rm (6)} are
satisfied for every $d\in C_G(\gp{a,b})$ of order $4$.

Moreover, the conditions {\rm 1)--5)} are also sufficient for
$K_{\lambda}G$  to be $f$-normal. The condition {\rm 6)} is
sufficient if $K$ is an integral domain
of characteristic $2$.
\endproclaim

\head {2. Lemmas} \endhead

Let $C_4$, $Q_8$ and $D_8$ be a cyclic group of order $4$,
a quaternion group of order $8$ and a dihedral group of
order $8$, respectively.  As usual, $x^y=y\m1xy$, $\exp(G)$
and   $C_G(\gp{a,b})$ denote the exponent of $G$ and the
centralizer of the subgroup $\gp{a,b}$ in $G$, respectively. 

It is easy to see that $\lambda_{g, g\m1}=\lambda_{g\m1, g}$
and $u_g\m1=\lambda_{g, g\m1}\m1u_{g\m1}$ hold for all $g\in G$.

\proclaim {Lemma 1}
The map $x \to x^f$ is  an involution of the ring $K_\lambda G$ if and only if
$$
f(gh)\lambda_{g, h}^2=f(g)f(h)\qquad   \text{for all}\quad  g, h\in G.
$$

\endproclaim

\demo {Proof}
Let the map $x \to x^f$ be  an involution of the ring $K_\lambda
G$. If $g, h \in G$, then $(u_gu_h)^f=u_h^fu_g^f$. Thus
$$
\align
\lambda_{g, h}f(gh)u_{gh}\m1  & =(\lambda_{g, h}u_{gh})^f=(u_gu_h)^f
=f(g)f(h)u_h\m1u_g\m1\\
& = f(g)f(h)(\lambda_{g, h}\m1u_{gh})\m1 \endalign $$
and  $f(gh)\lambda_{g, h}^2=f(g)f(h)$ for all $g, h\in G$.
\hfill$\qed$
\enddemo

Clearly, if  $K_\lambda G$ is a group ring, then the map $x
\to x^f$ is  an involution of the group ring $KG$  if and only if $f$
is a homomorphism from $G$ to $U(K)$.

\proclaim {Lemma 2}
If the ring $K_\lambda G$ is $f$-normal then the group $G$
satisfies one of the conditions  {\rm 1)--6)} of Theorem 1.
\endproclaim

\demo {Proof}
Let $K_\lambda G$ be an $f$-normal twisted group ring. If $a, b\in G$ and
$x=u_a+u_b\in K_\lambda G$, then $x^f=f(a)u_a\m1+f(b)u_b\m1$ and by (2)
$$
\gathered
f(a)\lambda_{a, a\m1}\m1\lambda_{a\m1, b}u_{a\m1b}+
f(b)\lambda_{b, b\m1}\m1\lambda_{b\m1, a}u_{b\m1a}\\  \vspace {3pt}
=f(a)\lambda_{a, a\m1}\m1\lambda_{b, a\m1}u_{ba\m1}+
f(b)\lambda_{b, b\m1}\m1\lambda_{a, b\m1}u_{ab\m1}. \endgathered \tag 7 $$
Now put $y=u_a(u_1+u_b)$. Then $y^f=(u_1+f(b)u_b\m1)f(a)u_a\m1$ and by (2)
$$
\lambda_{a, b}u_{ab}+f(b)\lambda_{b, b\m1}\m1\lambda_{a, b\m1}u_{ab\m1}=
\lambda_{b, a}u_{ba}+f(b)\lambda_{b, b\m1}\m1\lambda_{b\m1, a}u_{b\m1a}. \tag 8
$$
We shall treat two cases.
\smallskip
\item {I.}
Let  $[a, b]\ne 1$ for $a, b\in G$ and $a^2\ne 1$, $b^2\ne 1$. Then by (8)
$b^a=b\m1$ and by (7) $a^2=b^2$.  The factor system  satisfies
$$
\cases
\lambda_{a, b}=f(a)\lambda_{a, a\m1}\m1\lambda_{b, a\m1}=f(b)\lambda_{b,
b\m1}\m1\lambda_{b\m1, a}; \\  \vspace {3pt}
\lambda_{b, a}=f(a)\lambda_{a, a\m1}\m1\lambda_{a\m1, b}=f(b)\lambda_{b,
b\m1}\m1\lambda_{a, b\m1}.
\endcases \tag 9 $$

\item {II.}
Let $[a, b]\ne 1$ for $a, b\in G$ and $a^2=1$, $b^2\ne 1$. Then by (8) we have
$b^a=b\m1$ and by (7),
$f(a)=-\lambda_{a, a}$.  The factor system satisfies
$$
\cases
\lambda_{a, b}=f(b)\lambda_{b, b\m1}\m1\lambda_{b\m1, a}; \\ \vspace {3pt}
\lambda_{b, a}=f(b)\lambda_{b, b\m1}\m1\lambda_{a, b\m1}.
\endcases 
$$
Let $G$ be a nonabelian group and let $W=\{g\in G\mid g^2 \ne 1\}$.

First we consider the case when the elements of $W$  commute.
Then $\gp{w\mid w\in W}$ is an abelian subgroup and if $b\in W$ and
$a\in G\setminus \gp{W}$ then $a^2{=}1$ and  $(ab)^2=1$. Therefore,
$b^a=b\m1$ for all $b\in W$. Let $c\in C_G(\gp{W})\setminus\gp{W}$. Then
$c^2=1$, $(cb)^2=1$ and $cb\notin \gp{W}$.  But
$(cb)^2=c^2b^2=1$ and $b^2=1$, which is impossible. Therefore,
$C_G(\gp{W})=\gp{W}$ and $H=\gp{W}$ is a subgroup of index $2$. This
implies that $G= H\rtimes \langle a\rangle$ and $h^a=h\m1$ for all $h\in H$.

Now suppose that in $W$ there exist elements $a, b$ such that $[a,b]\ne 1$.
Since $a^2\ne 1$ and $b^2\ne 1$, by (I) we have
$a^2=b^2$ and $b^a=b\m1$. Then $b^2=ab^2a\m1=b^{-2}$ and the elements
$a, b$ are of order $4$. Clearly, the  subgroup $\gp{a, b}$ is a
quaternion  group of order $8$.  Let $c\in  C_G(\gp{a, b})$. If
$c^2\ne 1$ and $(ac)^2\ne 1$ then  (I) implies that
$(ac)^b=(ac)\m1$ and $c^2=1$, which is impossible. Therefore, if
$c\in C_G(\gp{a, b})$ then either $c^2=1$ or $c^2=a^2$.

Let $Q=\gp{a, b}$ be a quaternion  subgroup of order $8$ of $G$.
Then we will prove that $G=Q\cdot C_G(Q)$.
Suppose $g\in G\setminus C_G(Q)$. Pick the elements $a, b\in Q$
of order $4$ such that $a^g=a\m1$ and $b^g=b\m1$. Then $(ab)^g=ab$
and $d=gab\in C_G(Q)$. It  follows that $g=d(ab)\m1$ and
$G=Q\cdot C_G(Q)$. Similary as in \cite{3} we obtain that $G$ satisfies
the conditions 4) or 5) of the Theorem.
\hfill$\qed$
\enddemo

\head {3. Proof of Theorem} \endhead

{\it Necessity.\/}
Let $K_{\lambda }G$ be $f$-normal. Then by Lemma $2$ $G$ satisfies
one of the  conditions  1)--5) of the Theorem.

First, suppose that $G$ is abelian of exponent greater than $2$ and
$a,b\in G$. If $b^2\ne 1$ then by (8) we have $\lambda_{a, b}=\lambda_{b, a}$.

Let $a$, $b$ be  elements of order two and assume that there
exists  $c$ with $c^2=a$. Then by (1) we have
$$
\lambda_{c^2, b}\lambda_{c, c}=\lambda_{c, cb}\lambda_{c, b} \quad \text{and}
\quad \lambda_{b, c^2}\lambda_{c, c}=\lambda_{bc, c}\lambda_{b,c}. \tag 10 $$
Since $c^2\ne 1$, we have $\lambda_{c, cb}=\lambda_{bc, c}$ and
$\lambda_{c, b}=\lambda_{b, c}$. Then (10) implies
$\lambda_{c^2, b}=\lambda_{b, c^2}$  and $\lambda_{a, b}=\lambda_{b, a}$.

Let $a^2=b^2=1$ such that neither $a$ nor $b$ is the square of any
element of $G$. Then there exists $c$ such that
$(ca)^2\ne 1$. Thus,
$$
\lambda_{ca, b}\lambda_{c, a}=\lambda_{c, ab}\lambda_{a, b}, \quad
\lambda_{b, ac}\lambda_{a, c}=\lambda_{ba, c}\lambda_{b, a}.\tag 11 $$
Since $\lambda_{b, ac}=\lambda_{ac, b}$ and
$\lambda_{c,a}=\lambda_{a, c}$ from (11) we have $\lambda_{a, b}=
\lambda_{b, a}$ for all $a, b\in G$.  Therefore, if $G$ is abelian
and $G^2\ne 1$ then the factor system is symmetric and $K_\lambda
G$ is  commutative.

Now, let $\exp(G)=2$. Then by (8)
$\lambda_{a,b}+f(b)\lambda_{b,b}\m1\lambda_{a,b}
=\lambda_{b, a}+f(b)\lambda_{b, b}\m1\lambda_{b, a}$ for all $a,b\in G$.
Therefore, $(\lambda_{a,b}-\lambda_{b, a})(1+f(b)\lambda_{b, b}\m1)=0$.

Next, let $G=H\rtimes C_2$ be a semidirect product of an abelian group $H$
with $\exp(H)\ne 2$  and $C_2=\gp{a\mid a^2=1}$, and with $h^a=h\m1$ for all
$h\in H$.  Clearly, $K_\lambda H$ is  $f$-normal and the factor system of $H$
is symmetric. Put $x=u_h+u_a$ for  $h\in H$. Since $K_\lambda G$ is
$f$-normal, we have $S_f(x)=xx^f-x^fx=0$ and
$$
\gathered
f(a)\lambda_{a, a}\m1\lambda_{h, a}u_{ha}+ f(h)\lambda_{h,
h\m1}\m1\lambda_{a, h\m1}u_{ah\m1}\qquad \\ \vspace {3pt}
\quad\qquad -f(h)\lambda_{h,h\m1}\m1\lambda_{h\m1, a}u_{h\m1a}- f(a)\lambda_{a,a}\m1
\lambda_{a, h}u_{ah}=0.\endgathered \tag 12 $$
We will prove $u_au_h=u_h^fu_a$ for every $h\in H$.

First, let $h^2\ne 1$. Because  $h^a=h\m1$,   by (12)  we have
$$
u_a^fu_h+u_h^fu_a=0\tag  13  $$
and
$$
\cases
f(a)\lambda_{a, a}\m1\lambda_{a, h}+
f(h)\lambda_{h, h\m1}\m1\lambda_{h\m1, a}=0; \\ \vspace {3pt}
f(a)\lambda_{a, a}\m1\lambda_{h, a}+
f(h)\lambda_{h, h\m1}\m1\lambda_{a, h\m1}=0.
\endcases  \tag  14  
$$
Now, let $h^2=1$. Then there exists $b\in H$ with $b^2\ne 1$ and
$(hb)^2\ne 1$. Put  $x=u_a+u_hu_b$. Because $(hb)^a=(hb)\m1$ and
$S_f(x)=xx^f-x^fx=0$ we have
$$
u_a^fu_hu_b+(u_hu_b)^fu_a=0. \tag  15  $$
Since $[u_h, u_b]=1$,  by (15) and (13) we have
$u_a^f(u_hu_b)=u_a^fu_bu_h=-u_b^fu_au_h$ and
$u_a^f(u_hu_b)=-(u_hu_b)^fu_a=-u_b^fu_h^fu_a$.
Therefore, $u_au_h=u_h^fu_a$ for all $h\in H$ and this  implies
$$
\cases
\lambda_{a, h}=f(h)\lambda_{h, h\m1}\m1\lambda_{h\m1, a}; \\ \vspace {3pt}
\lambda_{h, a}=f(h)\lambda_{h, h\m1}\m1\lambda_{a, h\m1},
\endcases $$
and, by (14), $f(a)=-\lambda_{a, a}$.

Let $G$ be a hamiltonian $2$-group. It is well known (see \cite{5},
Theorem~12.5.4) that $G=Q_8\times W$, where $Q_8$ is a quaternion
group  and $\exp(W)|2$.
If $a, b\in G$ are noncommuting elements of order $4$, then
$a^b=a\m1$ and by (8) we have  4.i) of the theorem.
If $c, d\in G$ are involutions, then $c$ and $d$ commute with all
$a\in G$ of order $4$. Then $H=\gp{a, d, c}$ is abelian of exponent
greater than $2$ and $K_\lambda H$ is  $f$-normal. By the condition 1) of
the theorem, the factor system of $H$ is symmetric, and $u_a$ and
$u_b$ commute with $u_c$.

Now prove $f(c)=\lambda_{c, c}$ for all involutions $c\in G$.  Choose
the elements $a$, $b$ of order $4$ such that $b^a=b\m1$. Put $x=u_cu_a+u_b$.
Since $\lambda_{a, c}=\lambda_{c, a}$ and $\lambda_{b, c}=\lambda_{c, b}$ by
(2), for $x$ we obtain
$$
\align
S_f(x) & =(f(b)u_au_b\m1+f(a)f(c)\lambda_{c, c}\m1u_bu_a\m1\\ \vspace {3pt}
& \quad-f(b)u_b\m1u_a-f(a)f(c)\lambda_{c, c}\m1u_a\m1u_b)u_c=0 \endalign  $$
and $f(b)\lambda_{b, b\m1}\m1\lambda_{a, b\m1}=
f(c)f(a)\lambda_{c, c}\m1\lambda_{a, a\m1}\m1\lambda_{a\m1, b}$.
From this property and (9) we deduce  $f(c)=\lambda_{c, c}$.


Now, suppose that either $G=E\times W$ or $G=(E \Y C_4)\times W$, where
$E$ is an extraspecial $2$-group,  $\exp(W)|2$ and $E \Y C_4$  is
the central product of $E$  and $C_4=\gp{c}$ with $E'=\gp{c^2}$.

Let $a$ be an element of order $4$ and $h\in C_G(\gp{a})$.  Then
by the condition~1) of the theorem $\lambda_{a, h}=\lambda_{h, a}$.

Let   $\gp{a, b\mid a,b\in G}$ be  the quaternion subgroup   of order $8$.
Then by~4) we obtain (5).

Now, let $G=\gp{a,b}\Y \gp{d\mid d^4=1}$ be  a subgroup of $G$ and
$d^2=a^2$. Then $a^b=a\m1$,  and  $\gp{a, d}$ and $\gp{b, d}$ are
abelian subgroups of exponent not equal to $2$ and  by the condition 1) of the
theorem, $\lambda_{a, d}=\lambda_{d, a}$ and $\lambda_{b, d}=\lambda_{d, b}$.
Put $x=u_b+u_au_d$. Since $K_\lambda G$ is $f$-normal, we obtain
$$
\gather
f(b)\lambda_{b, b\m1}\m1\lambda_{a, b\m1}u_{ab\m1}u_d+
f(d)f(a)\lambda_{a, a\m1}\m1\lambda_{d, d\m1}\m1
\lambda_{b, a\m1}u_{ba\m1}u_{d\m1}\\ \vspace {3pt}
=f(d)f(a)\lambda_{a, a\m1}\m1\lambda_{d, d\m1}\m1
\lambda_{a\m1, b}u_{a\m1b}u_{d\m1}+
f(b)\lambda_{b, b\m1}\m1\lambda_{b\m1, a}u_{b\m1a}u_d\endgather $$
and by (5)
$$
\gather
\lambda_{b, a}\lambda_{ab\m1, d}u_{ab\m1d}+
f(d)\lambda_{d, d\m1}\m1\lambda_{a, b}\lambda_{ba\m1, d\m1}u_{ba\m1d\m1} \\
 \vspace {3pt}
=f(d)\lambda_{d, d\m1}\m1\lambda_{b, a}
\lambda_{a\m1b, d\m1}u_{a\m1bd\m1}+
\lambda_{a, b}\lambda_{b\m1a, d}u_{b\m1ad}.\endgather $$
Since $d^2\in G'$ and $a^2=b^2$, we have $a\m1bd\m1=abd$,
$ab\m1d=ba\m1d\m1$ and
$$
\lambda_{b, a}\lambda_{ba, d}+
f(d)\lambda_{d, d\m1}\m1\lambda_{a, b}\lambda_{ab, d\m1}=0. $$
Therefore, we proved  6.i).

If $\gp{a,b\mid a^4=b^2=1}$ is the dihedral subgroup of order
$8$ of $G$, then  by 3) of the theorem we have (4) and
$f(b)=-\lambda_{b,b}$.

Let  $L=D_8\Y C_4=\gp{a, b \mid a^4=b^2=1}\Y\gp{c}$.
Then  any  $x\in K_\lambda L$ can be written as
$x=x_0+x_1u_c$, where $x_0,x_1\in K_\lambda D_8$. Since $K_\lambda
G$ is $f$-normal,
$K_\lambda L$ is $f$-normal, too, and $(x_0^fx_1-x_1x_0^f)u_c=
(x_0x_1^f-x_1^fx_0)u_c^f$.  By the $f$-normality of $K_\lambda D_8$
$(x_0+x_1)(x_0+x_1)^f=(x_0+x_1)^f(x_0+x_1)$ and we have
$$
(x_0^fx_1-x_1x_0^f)u_c-
(x_0x_1^f-x_1^fx_0)u_c^f=(x_0^fx_1-x_1x_0^f)(u_c-u_c^f). $$
If $x_0^fx_1-x_1x_0^f$ can be written as a sum of elements of form
$u_a^fu_b-u_bu_a^f$ then
$$
\gather
(x_0^fx_1-x_1x_0^f)(u_c-u_c^f)=(\lambda_{b,a}\lambda_{ba,c}
+f(c)\lambda_{c,c\m1}\m1\lambda_{a,b}\lambda_{ab,c\m1})u_{bac}\\
-(\lambda_{a,b}\lambda_{ab,c}+
f(c)\lambda_{c,c\m1}\m1\lambda_{b,a}\lambda_{ba,c\m1})u_{abc}=0\endgather $$
and we have (6).

\subhead
Sufficiency
\endsubhead
We wish to prove that $S_f(x)=xx^f-x^fx$ is equal to $0$ for all $x\in KG$.
Let $x=\sum_{g\in G}\alpha_gu_g\in K_\lambda G$.
It is easy to see that $S_f(x)$ is a sum of elements of the  form
$$
\align
S_f(g,h) & =\alpha_g\alpha_h(f(h)\lambda_{h, h\m1}\m1
 \lambda_{g,h\m1}u_{gh\m1}+f(g) \lambda_{g, g\m1}\m1
 \lambda_{h, g\m1}u_{hg\m1}\\  \vspace {3pt}
&\quad-f(h)\lambda_{h, h\m1}\m1\lambda_{h\m1, g}u_{h\m1g}-
f(g) \lambda_{g, g\m1}\m1\lambda_{g\m1, h}u_{g\m1h}). \endalign $$
First, let $G$ be abelian of exponent greater than $2$, and assume that the
factor system of $G$ is symmetric. Then $K_{\lambda}G$ is commutative, and
therefore, $f$-normal.

Next, suppose that $G$ is of exponent $2$  and the factor system satisfies
$$
(\lambda_{g, h}-\lambda_{h, g})(1+f(h)\lambda_{h, h}^{-1})=0\quad 
\text{ for all}\quad  g, h\in G. 
$$
This  implies
$(\lambda_{g,h}-\lambda_{b,h})(f(g)\lambda_{g,g}\m1-f(h)\lambda_{h, h}\m1)=0$
for all $g, h\in G$.  Then
$$
\gather
S_f(g, h)=\alpha_g\alpha_h( f(h)\lambda_{h, h}\m1\lambda_{g,
h}u_{gh}+ f(g)\lambda_{g, g}\m1\lambda_{h, g}u_{hg}-
f(h)\lambda_{h, h}\m1\lambda_{h, g}u_{hg}\\  \vspace {3pt}
-f(g)\lambda_{g, g}\m1\lambda_{g, h}u_{gh})=
\alpha_g\alpha_h (f(h)\lambda_{h, h}\m1- f(g)\lambda_{g, g}\m1)
(\lambda_{g, h}-\lambda_{h, g})u_{gh}=0 \endgather $$
and $S_f(x)=0$, thus,  $K_\lambda G$ is $f$-normal.

Now, let $G=H\rtimes C_2$, where $H$ is an abelian group of
exponent not equal to $2$  and $C_2=\gp{a}$ with $h^a=h\m1$
for all $h\in H$.  Using the properties  of the factor system
we obtain
$$
\gathered
f(a)u_a\m1u_h=-f(h)u_h\m1u_a,\qquad f(a)u_hu_a\m1=-f(h)u_au_h\m1,\\
 \vspace {3pt}
u_a^fy=-y^fu_a,\qquad yu_a^f=-u_ay^f \endgathered  \tag 16$$
for any $h\in H$ and $y\in K_\lambda H$.
If $x=x_1+x_2u_a\in K_\lambda G\,$ where $x_1, x_2\in K_\lambda H$,
then $x^f=x_1^f+f(a)u_a\m1x_2^f$ and
$$
xx^f=x_1x_1^f+f(a)x_1u_a\m1x_2^f +x_2u_ax_1^f+f(a)x_2x_2^f. $$
Because in $K_\lambda H$ the factor
system is symmetric and $K_\lambda H$ is commutative, by (16) we have
$$
xx^f=x_1x_1^f+(x_2x_1-x_1x_2)u_a+f(a)x_2x_2^f=x_1x_1^f+f(a)x_2x_2^f. $$
Similarly, $x^fx=x_1^fx_1+f(a)x_2^fx_2$ and we conclude that $S_f(x)=0$
and $K_\lambda G$ is $f$-normal.

Next, let  $G$ be a hamiltonian $2$-group.  Then $G=Q_8\times W$, where
$Q_8=\gp{a, b}$ is  a quaternion group and $\exp(W)|2$. Suppose  that the
conditions 4.i)--4.ii) of the theorem are satisfied.  If $H=\gp{a^2, W}$
then any element $x\in K_\lambda G$ can be written as
$$
x=x_0+x_1u_a+x_2u_b+x_3u_{ab},  $$
where $x_i\in K_\lambda H$, ($i=0,\dots, 3$).  Since $\gp{a}\times H$
and $\gp{b}\times H$ are  abelian groups of exponent $4$, by the
condition 1) of the theorem the elements $x_0$, $x_1$, $x_2$, $x_3$
commute with $u_a$, $u_b$ and $u_{ab}$.  Since $K_\lambda H$ is
$f$-normal, we have
$x_ix_j^f-x_i^fx_j=x_j^fx_i-x_jx_i^f$.  Using these properties we obtain
$$ \spreadlines{3pt}
\align
S_f(x) & =(x_1x_2^f-x_1^fx_2)(\lambda_{b, a}u_{ba}-\lambda_{a,
b}u_{ab})\\
&\quad + (x_1x_3^f-x_1^fx_3) (\lambda_{ab, a}u_{b}-\lambda_{a,ab}u_{b^3})\\
& \quad+(x_2x_3^f-x_2^fx_3) (\lambda_{ab,b}u_{a^3}-\lambda_{b, ab}u_{a}).
\endalign $$

Clearly, the element $x_ix_j^f- x_i^fx_j$
can be written as a sum of elements of form
$$
S_f(c, d)=\gamma_{c, d}(f(d)u_cu_d\m1-f(c)u_c\m1u_d), $$
where $c, d\in H$.
Since $H$ is an elementary $2$-subgroup, by the condition 4.ii)
$f(d)=\lambda_{d, d}$, $f(c)=\lambda_{c, c}$, and we obtain
$$
S_f(c, d)=\gamma_{c, d}(f(d)\lambda_{d, d}\m1
\lambda_{c, d}u_{cd}-f(c)\lambda_{c, c}\m1\lambda_{c, d}u_{cd})=0. $$
Therefore, $S_f(x)=0$ and $K_\lambda G$ is  $f$-normal.

Next, let $G=H\times W$, where $H$ is an extraspecial $2$-group and
$\exp(W)|2$. Since $G$ is a locally finite group, it  suffices to
establish the  $f$-normality of all finite subgroups $H$ of $G$.
Let $G$ be a finite group and $G=H\times W$, where $H$ is a finite
extraspecial $2$-group and $\exp(W)|2$.  We know (see \cite{4},
Theorem III.13.8) that  $H$ is a central product of $n$ copies of
dihedral groups of order $8$ or a central product of a  quaternion
group of order $8$ and $n-1$  copies of  dihedral groups of order
$8$. We can write $H_n=H$. Then $G=H_n\times W$ and by induction on
$n$ we prove the  $f$-normality of $K_\lambda G$.

If $n=1$ then either $H_1=Q_8$ or $H_1=D_8$ or $H_1=Q_8\Y C_4$.
In the first and second cases  the $f$-normality $K_\lambda G$ is
implied by the conditions 3) or 4)  of the theorem.

Let $G=Q_8\Y C_4$. Then any element $x\in K_\lambda G$ can be
written as $x=x_0+x_1u_c$, where $x_i\in K_\lambda Q_8$, $c\in C_4$
and $c^2\in Q_8$. From the $f$-normality of $K_\lambda Q_8$ we obtain
$x_0^fx_1-x_1x_0^f=x_1^fx_0-x_0x_1^f$ and
$S_f(x)=(x_0^fx_1-x_1x_0^f)(u_c-u_c^f)$. The element $x_0^fx_1-x_1x_0^f$
can be written as a sum of elements of form $\alpha(u_a^fu_b-u_bu_a^f)$,
where $\alpha\in K$, $a,b\in Q_8$.
We will prove $S_f(a,b)=(u_a^fu_b-u_bu_a^f)(u_c-u_c^f)=0$ for all $a,b\in Q_8$.

If $a,b\in Q_8$ does  not generate  $Q_8$ then
$u_au_b=u_bu_a$ and $S_f(a,b)=0$. Let $\gp{a,b}=Q_8$. Then by (5)
$$\spreadlines{3pt}
\align
S_f(a,b) & =(\lambda_{b,a}u_{ba}-\lambda_{a,b}u_{ab})(u_c-u_c^f)\\
& =(\lambda_{b,a}\lambda_{ba,c}+f(c)\lambda_{c,c\m1}\m1 \lambda_{a,b}
\lambda_{ab,c\m1})u_{bac}\\
&\quad +(\lambda_{a,b}\lambda_{ab,c}+f(c)\lambda_{c,c\m1}\m1 \lambda_{b,a}
\lambda_{ba,c\m1})u_{abc} \endalign $$
and from (6) $S_f(a,b)=0$.

It is easy to see $D_8 \Y D_8\cong   Q_8 \Y  Q_8$,
and $H_n$ $(n>1)$  can be written as $Q_8 \Y H_{n-1}$.

Let $Q_8=\gp{a, b}$ and $L=W\times H_{n-1}$. Any element
$x\in K_\lambda G$ can be written as
$$
x=x_0+x_1u_a+x_2u_b+x_3u_au_b,  $$
where $x_i\in K_\lambda L$. By 6.i) the   $x_i$  commute with
$u_a$ and $u_b$. Since
$\gp{a, b}$ is a quaternion group of order $8$, by the condition
\rm{ 6.ii)} of the theorem we have $u_au_b=u_b^fu_a=u_bu_a^f$.
Hence,
$$\spreadlines{3pt}
\alignat 1
S_f(x)  ={}&(x_0x_1^f{-}x_1^fx_0)u_a^f{+}(x_0x_2^f{-}x_2^fx_0)u_b^f
{+}(x_0x_3^f{-}x_3^fx_0)u_b^fu_a^f\\
&\ {+}(x_1x_0^f{-}x_0^fx_1)u_a{+}(x_1x_2^f{-}x_1^fx_2)u_au_b^f
 {+}(x_1x_3^f{-}x_1^fx_3)u_bf(a)\\
&\ {+}(x_2x_0^f{-}x_0^fx_2)u_b{+}(x_2x_1^f{-}x_2^fx_1)u_au_b{+}(x_2x_3^f
{-}x_2^fx_3)u_a^ff(b) \tag 17 \\
&\ {+}(x_3x_0^f-x_0^fx_3)u_au_b+
(x_3x_1^f-x_3^fx_1)u_au_{ab}\\
&\  {+}(x_3x_2^f{-}x_3^fx_2)u_af(b). \endalignat $$

Since by induction $K_\lambda L$ is $f$-normal, $(x_i+x_j)(x_i+x_j)^f=$
\newline$(x_i+x_j)^f(x_i+x_j)$ implies
$x_ix_j^f-x_i^fx_j=x_j^fx_i-x_jx_i^f$ and
$x_ix_j^f-x_j^fx_i=x_i^fx_j-x_jx_i^f$.
Therefore,  by (17)
$$\spreadlines{3pt}
\align
S_f(x) & =(x_0x_1^f-x_1^fx_0)(u_a^f-u_a)+(x_0x_2^f-x_2^fx_0)(u_b^f-u_b)\\
&\quad +(x_0x_3^f-x_3^fx_0)(u_{a}^f-u_a)u_b+(x_1x_2^f-x_1^fx_2)u_a(u_b^f-u_b)\\
& \quad +(x_1x_3^f-x_1^fx_3)u_a(u_b-u_b^f)f(a)+(x_2x_3^f-x_2^fx_3)
(u_a^f-u_a)f(b). \endalign $$


Clearly, the element $x_ix_j^f- x_j^fx_i$
can be written as a sum of  elements of form
$S_f(c, d)=\gamma_{c, d}(u_cu_d^f-u_d^fu_c)$,
where $c, d\in L$, $\gamma_{c, d}\in K$. We will  prove
$S_f(c,d,a)=(u_cu_d^f-u_d^fu_c)(u_a-u_a^f)=0$ for any  $c, d\in L$.

\medskip
We consider the following cases:

\smallskip
{\it Case 1).}
Let $[c,d]=1$. Then $L=\gp{c,d,a}$ is abelian with $\exp(L)\ne 2$, and by 6.i)
 the factor system is symmetric and $S_f(c,d,a)=0$.

\smallskip
{\it Case 2).}
Let $\gp{c,d}=Q_8$. Then by 6.ii) (5) holds  and
$$
\split
(u_cu_d^f-u_d^fu_c)(u_a-u_a^f)=&(\lambda_{d,c}\lambda_{dc,a}+
f(a)\lambda_{a,a\m1}\m1\lambda_{c,d}\lambda_{cd,a\m1})u_{dca}\\
-&(\lambda_{c,d}\lambda_{cd,a}+f(a)\lambda_{a,a\m1}\m1\lambda_{d,c}
\lambda_{a\m1,dc})u_{cda}
\endsplit
$$

Now by 6.ii) the property  (6) is satisfied and we conclude $S_f(c,d,a)=0$.

\smallskip
{\it Case 3).}
Let $\gp{c,d}=D_8$ and  $c^4=d^2=1$. Then by 6.iii) $f(d)=-\lambda_{d,d}$ and
by (4) we have that 
$$\spreadlines{3pt}
\align
(u_cu_d^f & -u_d^fu_c)(u_a-u_a^f)=
(\lambda_{d,c}u_{dc}-\lambda_{c,d}u_{cd})(u_a-u_a^f)\\
& =(\lambda_{c,d}\lambda_{cd,a}+
f(a)\lambda_{a,a\m1}\m1\lambda_{dc,a\m1}\lambda_{d,c})u_{cda}\\
& \quad+(\lambda_{d,c}\lambda_{dc,a}+
f(a)\lambda_{a,a\m1}\m1\lambda_{cd,a\m1}\lambda_{c,d})u_{dca}.\endalign  $$
Now by 6.ii) we have  (6) and we conclude $S_f(c,d,a)=0$.

\smallskip
{\it Case 4).}
Let $\gp{c,d}=D_8$ and  $d^4=c^2=1$. Then  by (4)
$$
\align
u_cu_d^f-u_d^fu_c & = f(d)\lambda_{d,d\m1}\m1 \lambda_{c,d\m1}u_{dc}-
f(d)\lambda_{d,d\m1}\m1 \lambda_{d\m1,c}u_{cd}\\
& = \lambda_{d,c}u_{dc}-\lambda_{c,d}u_{cd}.\endalign  $$
Similarly to the case 3) we have $S_f(c,d,a)=0$.

\smallskip
{\it Case 5).}
Let $\gp{c,d}=D_8$ and  $d^2=c^2=1$. Then by 6.iii) $f(d)=-\lambda_{d,d}$.
In $\gp{c,d}$ we  choose a  new generator system
$\{a_1,b_1\mid a_1^4=b_1^2=1, a_1^{b_1}=a_1\m1\}$ such that
$c=b_1$ and  $d=a_1^ib_1$, where $i=1$ or $3$. Then $a^2=a_1^2$ and
$$
\gather
(u_cu_d^f-u_d^fu_c)(u_a-u_a^f)=(u_du_c-u_cu_d)(u_a-u_a^f)\\
=\lambda_{a_1^i,b_1}\m1 (u_{a_1^i}u_{b_1}-u_{b_1}u_{a_1^i})
(u_a-u_a^f)u_{b_1}.\endgather $$
As in the Case 3) it is easy to see $(u_{a_1^i}u_{b_1}-u_{b_1}u_{a_1^i})
(u_a-u_a^f)=0$ and  $S_f(c,d,a)=0$.

Analogously, the element $x_ix_j^f- x_i^fx_j$ can be written as a sum of
elements of form $\gamma_{c, d}(u_cu_d^f-u_c^fu_d)$, where $c,d\in L$. Let us
prove that if $c, d\in L$, then 
$$
S_f(c,d,a)=(u_cu_d^f-u_c^fu_d)(u_a-u_a^f)=0.
$$

Let $z\in L$,  $a\in Q_8$ be commuting elements of  order $4$ with $z^2=a^2$.
First, we will  prove that $K$ is of characteristic $2$,  then
$(u_z+u_z^f)(u_a+u_a^f){=}0$.

Indeed,
$$
\gather
(u_z+u_z^f)(u_a+u_a^f)=(\lambda_{z,a}+f(z)\lambda_{z,z\m1}
\m1f(a)\lambda_{a,a\m1}\m1\lambda_{z\m1,a\m1})u_{za}\\ \vspace {3pt}
+(f(a)\lambda_{a,a\m1}\m1\lambda_{z,a\m1}+
f(z)\lambda_{z,z\m1}\m1\lambda_{z\m1,a})u_{za^3}.\endgather $$
First let  $za$ be  a noncentral element of order $2$. Then by 6.iii)
$f(za)=\lambda_{za,za}$. Since $((u_zu_a)u_a)u_{a^3}=u_z(u_a(u_au_{a^3}))$ we
conclude that
$$
\lambda_{z,a}\lambda_{za,a}\lambda_{za^2,a^3}=
\lambda_{z,a}\lambda_{a,1}\lambda_{a,a\m1}$$
and $\lambda_{a,a\m1}\m1=\lambda_{z^3,a^3}\m1\lambda_{za,a}\m1$.
Clearly,
$f(z)f(a)=f(za)\lambda_{z,a}^2=\lambda_{za,za}\lambda_{z,a}^2$ and
$$
\spreadlines {3pt}
\alignat 1
\lambda_{z,a}+ f(z) & \lambda_{z,z\m1}\m1f(a)\lambda_{a,a\m1}\m1
  \lambda_{z\m1,a\m1}\\
& =\lambda_{z,a}(1+(\lambda_{za,az}\lambda_{a,z})\lambda_{z,z\m1}\m1
 \lambda_{a,a\m1}\m1\lambda_{z\m1,a\m1})\\
& =\lambda_{z,a}(1+\lambda_{z,za^2}\lambda_{a,az}\lambda_{a,a\m1}\m1
\lambda_{z,za^2}\m1\lambda_{z\m1,a\m1} \tag 18 \\
& =\lambda_{z,a}(1+(\lambda_{a,az}\lambda_{zaa,a\m1})\lambda_{a,a\m1}\m1)\\
& =\lambda_{z,a}(1+
\lambda_{za,aa\m1}\lambda_{a,a\m1}\lambda_{a,a\m1}\m1)=2\lambda_{z,a}=0.
\endalignat $$

By (1) we have
$$
\gather
(\lambda_{z,a\m1}\lambda_{za\m1,za\m1})\lambda_{z\m1,a}=
\lambda_{z,a\m1za\m1}\lambda_{a\m1,z\m1a}\lambda_{z\m1,a}\\ \vspace {2pt}
=\lambda_{z,z\m1}(\lambda_{a\m1,az\m1}\lambda_{a,z\m1})
=\lambda_{z\m1,z}\lambda_{aa\m1,z\m1}\lambda_{a,a\m1}=
\lambda_{z\m1,z}\lambda_{a,a\m1} \endgather $$
and since $az\m1$ has order $2$, $f(az\m1)=\lambda_{az\m1,az\m1}$,
and we obtain
$$
\spreadlines{2pt}
\alignat 1
f(a & \m1)\m1  f(a\m1)(f(a)\lambda_{a,a\m1}\m1\lambda_{z,a\m1}+
f(z)\lambda_{z,z\m1}\m1\lambda_{z\m1,a})\\
& =f(a\m1)\m1(\lambda_{a,a\m1}^2\lambda_{a,a\m1}\m1\lambda_{z,a\m1}
{+}f(az\m1)\lambda_{a\m1,z}^2\lambda_{z,z\m1}\m1\lambda_{z\m1,a})\\
& =f(\!a\m1\!)\!\m1\!(\!\lambda_{a,a\m1}\lambda_{z,a\m1}{+}
\lambda_{z,a\m1}\!(\!\lambda_{z,a\m1}
\lambda_{az\m1,az\m1}\lambda_{z\m1,a}\!)\!\lambda_{z,z\m1}\m1\!)\! \tag 19 \\
& =f(a\m1)\m1(\lambda_{z,a\m1}(\lambda_{a\m1,a}-\lambda_{z,z\m1}
\lambda_{a\m1,a}\lambda_{z\m1,z}\m1\lambda_{z\m1,a})\\
& =2f(a\m1)\m1\lambda_{z,a\m1}\lambda_{a,a\m1}=0.\endalignat $$
Clearly, if $[c,d]=1$ then $S_f(c,d,a)$ can be written as
$$
\gathered
S_f(c,d,a)=(u_cu_d^f+(u_d^fu_c)^f)(u_a-u_a^f)\\
=f(d)\lambda_{d,d\m1}\lambda_{c,d\m1}(u_{cd\m1}-u_{cd\m1}^f)(u_a-u_a^f).
\endgathered \tag 20 $$
Similarly,  the element $x_ix_j^f- x_i^fx_j$
can be written as a sum of elements of form
$\gamma_{c, d}(u_cu_d^f-u_c^fu_d)$,
where $c, d\in L$. Now let us prove
$$
S_f(c,d,a)=(u_cu_d^f-u_c^fu_d)(u_a-u_a^f)=0,
$$ 
where $c, d\in L$.

\medskip
We consider the following cases:

\smallskip
{\it Case 1).}
Let $[c, d]=1$, $c^2=d^2=1$ and $c,d\notin\zeta(G)$.  Then
$S=\gp{c,d,a}$ is abelian of exponent greater that $2$ and by 6.i) the factor
system of $S$ is symmetric.  We know that in $L$ every element of order
$2$ is either  central or coincides with  a noncentral element of some
dihedral subgroup of order $8$.  Since $c,d\notin\zeta(G)$, we
have $f(c)=\lambda_{c,c}$ and $f(d)=\lambda_{d,d}$ and
$$
S_f(c,d,a)=
\lambda_{c,d}(f(d)\lambda_{d,d}\m1-
f(c)\lambda_{c,c}\m1)u_{cd}(u_a-u_a^f)=0.$$

\smallskip
{\it Case 2).}
Let $[c, d]=1$, $c^2=d^2=1$ and $c,d\in\zeta(G)$. Then $c=d=a^2$ and
$S_f(c,d,a)=0$.

\smallskip
{\it Case 3).}
Let $[c, d]=1$, $c^2=d^2=1$ and $c\in \zeta(G)$,
$d\notin\zeta(G)$. Then $f(d)=\lambda_{d,d}\m1$, $c=a^2$ and
$$
\gather
S_f(c,d,a)=-u_d(u_{a^2}+u_{a^2}^f)(u_a-u_a^f)\\
=-u_d(\lambda_{a,a^2}u_{a\m1}-f(a)\lambda_{a,a\m1}\m1u_a)
(1+f(a^2)\lambda_{a^2,a^2}\m1). \endgather $$
Since $K$ is an integral domain of characteristic $2$ and
$f^2(a^2){=}\lambda_{a^2,a^2}^2f(a^4){=}\mathbreak\lambda_{a^2,a^2}^2$,  we
conclude $f(a^2)=\pm \lambda_{a^2,a^2}$ and $S_f(c,d,a)=0$.

\smallskip
{\it Case 4).}
Let $[c, d]=1$, $d^2=1$ and suppose that $c$ has order $4$.
Then $dc$ has order $4$ and by  (20) $S_f(c,d,a)=0$.

\smallskip
{\it Case 5).}
Let $[c, d]=1$ with  $c,d$ of order $4$. Then $d^2=c^2=a^2$,
$$
S_f(c,d,a)=(f(d)\lambda_{d,d\m1}\m1\lambda_{c,d\m1}
+f(c)\lambda_{c,c\m1}\m1\lambda_{c\m1,d})u_{cd\m1}(u_a-u_a^f), $$
and by (19) we have $S_f(c,d,a)=0$.

\smallskip
{\it Case 6).}
Let $\gp{c,d}$ be a quaternion group of order $8$. Then
by 6.ii) (5) holds and
$$
\align
u_cu_d^f-u_c^fu_d & =
(f(d)\lambda_{d,d\m1}\m1\lambda_{c,d\m1}-
f(c)\lambda_{c,c\m1}\m1\lambda_{c\m1,d})u_{c\m1d}\\
& =(\lambda_{d,c}-\lambda_{d,c})u_{c\m1d}=0. \endalign $$

\smallskip
{\it Case 7).}
Let $\gp{c,d}\cong D_8$. If
$c^2\ne 1$, then   $f(d)=\lambda_{d,d}$ and
$$
\gather
S_f(c,d,a)=(\lambda_{c,d}u_{cd}
+f(c)\lambda_{c,c\m1}\m1\lambda_{c\m1,d}u_{dc})(u_a-u_a^f)\\
=(\lambda_{c,d}u_{cd}+\lambda_{d,c}u_{dc})(u_a-u_a^f)
=(\lambda_{c,d}\lambda_{cd,a}+
f(a)\lambda_{a,a\m1}\lambda_{d,c}\lambda_{dc})u_{acd}\\
-(\lambda_{d,c}\lambda_{dc,a}+
f(a)\lambda_{a,a\m1}\lambda_{c,d}\lambda_{cd,a\m1})u_{adc}.\endgather $$

By (6) we obtain $S_f(c,d,a)=0$.

\smallskip
{\it Case 8).}
Let $\gp{c,d}$ be a dihedral group of order $8$ and
$c^2=d^2=1$. Then $f(d)=\lambda_{d,d}$, $f(c)=\lambda_{c,c}$ and
$S_f(c,d,a)=2u_cu_d(u_a-u_a^f)=0$.
\hfill$\qed$

\enddocument